\documentclass[a4paper,12pt,twoside]{article}

\usepackage{amssymb,epsfig}
\usepackage{amsmath}
\usepackage{amscd}

\usepackage{fancyheadings}
\newtheorem{thm}{Theorem}[section]
\newtheorem{lem}[thm]{Lemma}
\newtheorem{prop}[thm]{Proposition}
\newtheorem{cor}[thm]{Corollary}

\newtheorem{defn}[thm]{Definition}

\newtheorem{rmk}[thm]{Remark}

\newcommand{\qed}{\hfill $\Box$ \vspace{.5cm}}
\newcommand{\pf}{{\bf Proof. }}

\newcommand{\ind}{\operatorname{ind}}
\newcommand{\coker}{\operatorname{coker}}
\newcommand{\Hom}{\operatorname{Hom}}

\newcommand{\topind}{\operatorname{top-ind}}
\newcommand{\supp}{\operatorname{supp}}
\newcommand{\Img}{\operatorname{Im}}

\newcommand{\psos}{pseudodifferential operators }
\newcommand{\pso}{pseudodifferential operator }
\newcommand{\psimult}{\Psi^0_{mult}}
\newcommand{\spinc}{spin$^c$}

\newcommand{\into}{\longrightarrow}

\renewcommand{\hat}{\widehat}
\renewcommand{\tilde}{\widetilde}

\newcommand{\NN}{\mathbb{N}}
\newcommand{\ZZ}{\mathbb{Z}}

\newcommand{\RR}{\mathbb{R}}
\newcommand{\CC}{\mathbb{C}}

\newcommand{\K}{\mbox{$\cal K$}}

\graphicspath{{./}{Figs/}{Plots/}}

\title{A $K$-Theory Proof of the \\Cobordism Invariance
of the Index}
\author{Catarina Carvalho\thanks{Supported by Funda\c{c}\~{a}o para a Ci\^{e}ncia e Tecnologia, PRAXIS XXI/BD/19876/99, and by NWO through the Pioneer Project 616.062.384 of N. P. Landsman.}\\
Korteweg-de Vries Institute, University of Amsterdam\\ email: {\tt ccarvalh@science.uva.nl}}

\begin{document}

\maketitle

\begin{abstract}
We give a proof of the cobordism invariance of the index of  elliptic pseudodifferential operators  on  $\sigma$-compact manifolds, where, in the non-compact case, the operators are assumed to be multiplication outside a compact set. 
We show that, if the principal symbol class of such an elliptic operator on the boundary of a manifold $X$ has a suitable extension to $K^1(TX)$, then its index is zero. This condition is incorporated into the definition of a  cobordism group for non-compact manifolds, called here ``cobordism of symbols''.
 Our proof is topological, in that we use properties of the push-forward map in $K$-theory defined by Atiyah and Singer, to reduce it to $\RR^n$. 
In particular, we generalize  the invariance of the index with respect to the push-forward map to the non-compact case, and obtain an extension of the $K$-theoretical index formula of Atiyah and Singer to operators that are multiplication outside a compact set. Our results hold also for $G$-equivariant operators, where $G$ is a compact Lie group.
\end{abstract}

\section*{Introduction}\label{s int}

 Atiyah and Singer originally used the cobordism invariance of the index of twisted signature operators on closed manifolds as the main ingredient in their first proof of the index theorem \cite{a-s-cobd}. 
 Since then, many proofs of cobordism invariance have been given for Dirac operators on closed manifolds (see for instance \cite{braverman-cobdsm,higson-cobdsm}), using mainly geometrical arguments that rely on the specific structure of the Dirac operator. 

Here, we give a topological proof of the cobordism invariance of the index at the level of $K$-theory by establishing a sufficient condition on the principal symbol class.  One of the main points of this paper is that the cobordism invariance of the index is a consequence of the compatibility of the index with push-forward maps. This approach has the  advantage of applying to general elliptic \psos and, moreover, of holding also on non-compact manifolds, where we consider 
operators that are multiplication at infinity (as in Definition \ref{mpsos}). 

More precisely, for a $\sigma$-compact manifold $M$ without boundary, we first show in Section \ref{s an-ind} that we have a well-defined index map 
$$\ind: K^0(TM)\to\ZZ,$$
which computes the indices of elliptic \psos that are multiplication at infinity on $M$. If $M$ is the boundary of some $\sigma$-compact manifold $X$, we use Bott periodicity to define a map of restriction of symbols (see Section \ref{s cobd inv}),
$$u_M : K^1(TX) \rightarrow K^0(TM).$$
Our main result, Theorem \ref{cobdsm inv}, states that, if $\sigma\in \Img(u_M)$ then $\ind (\sigma)=0$.
 Note that, for any manifold $M$, there is always a (non-compact) manifold $X$ such that $\partial X=M$,
 and therefore in order to obtain a meaningful notion of non-compact cobordism one has to require that some extra structure is preserved.  
On the other hand, the condition that $\sigma\in \Img(u_M)$ is as much a condition on the symbol as on the manifold: we want $M$ to be the boundary of a manifold $X$ satisfying prescribed conditions. These facts motivate the incorporation of the $K$-theoretical condition given above into the definition of a non-compact cobordism, so-called \textit{cobordism of symbols}, on the set of pairs $(M,\sigma)$, where $M$ is a $\sigma$-compact manifold and $\sigma\in K^0(TM)$. We say that $(M,\sigma)$ is  \textit{symbol cobordant to zero} if there exist a manifold $X$ and $\omega\in K^1(TX)$ such that

(i) $\partial X= M$

(ii) $u_{M}(\omega)=\sigma$.

\noindent In the compact case, this relation  is closely related to the definition of cobordism for twisted signature operators defined by Atiyah and Singer in \cite{a-s-cobd} (see also \cite{palais}). The cobordism invariance of the index can now be stated as follows:
\newline

\textit{Let $P$ be an order zero bounded elliptic \pso that is multiplication at infinity on a $\sigma$-compact manifold $M$, with principal symbol $\sigma\in K^0(TM)$. If $(M,\sigma)$ is symbol cobordant to zero then $\ind(P)=0$.}
\newline

The main tool used in the proof of the above result is the push-forward map in $K$-theory, defined in \cite{a-sI}. To an embedding $i: X\to Y$\footnote{For manifolds with boundary, the embedding is required to be \textit{neat} (see Section \ref{s p-fwd map}).}, Atiyah and Singer associated a 'wrong-way functoriality' map
$$i_! : K^0(TX)\to K^0(TY),$$
which was used to give a purely $K$-theoretical proof of the index theorem; this relied heavily on the fact that, for closed manifolds, $\ind\circ i_!=\ind$.

 We prove that cobordism of symbols is preserved by the push-forward map, which reduces the proof to $\RR^n$; in this case, this is a simple computation in $K$-theory. We use crucial functoriality properties of the push-forward map, namely its invariance under excision, and,  more importantly, the fact that the index is invariant with respect to push-forward; these properties are given here in the more general setting of non-compact manifolds (possibly with boundary). In particular, one obtains a generalization of the $K$-theoretical formula for the index, given in \cite{a-sI}, to operators that are multiplication at infinity (see Section \ref{s ind form}).

For compact manifolds, an equivalent $K$-theoretical formulation of cobordism invariance is given in  \cite{moroianu-cobdsm}, where Moroianu gives also a condition for cobordism invariance at the level of operators. The methods of proof are however entirely different from ours, using the analysis of cusp \psos on manifolds with boundary and the calculus of non-commutative residues. 

To our knowledge, there is only one other notion of non-compact cobordism in the index theory setting. In \cite{braverman-cobdsm}, Braverman considered  a notion of non-compact cobordism, along the lines of \cite{GGK_book, karshon_ncpt_cobdsm}, for equivariant Dirac operators, and showed that the transverse index of such operators is invariant under this cobordism relation. 
%
 While it is not clear at this point what is the connection, if any, with the notion of non-compact cobordism presented here, our point of view seems to have the advantage of adapting better to our purposes.

We now review the contents of each section.  In Section \ref{s an-ind}, we define the class of \psos that are a multiplication at infinity (Definition \ref{mpsos}) and give its main properties. Namely, we give conditions for boundedness  and Fredholmness (Propositions \ref{psos bdd} and \ref{fred}) and show that there is indeed a well-defined index map in $K$-theory. We also show that this index map is invariant under excision (Proposition \ref{an-ind vs ext}).

  Our main result, Theorem \ref{cobdsm inv}, is stated in Section \ref{s cobd inv}. We first establish the notion of cobordism of symbols and show that the twisted signature operator on the boundary of a compact manifold is cobordant to zero. 

In Section \ref{s p-fwd map}, we consider the push-forward map in $K$-theory, defined here for manifolds with boundary. We give some functoriality properties, namely its behavior with respect to the Bott map, and an excision property (Propositions \ref{i_! vs ext} and \ref{i_! vs Bott}). We also prove functoriality with respect to restriction to the boundary, in Proposition \ref{i_! vs rest to bdry}, which is relevant in order to deal with cobordism invariance. To finish the section, we show that the crucial invariance of the index with respect to the push-forward map can be extended to our non-compact setting (Theorem \ref{i_! vs an-ind}).

In Section \ref{s pf cbd inv} we prove Theorem \ref{cobdsm inv}. At this point, one can easily show that the relation of cobordism of symbols is preserved by the push-forward map (Proposition \ref{i_! vs cobdsm symb}) and then it suffices to compute the index for operators on $\RR^n$.  
Finally, in Section \ref{s ind form}, we show how the invariance of the index with respect to the push-forward map easily  yields a $K$-theoretical index formula for operators that are multiplication at infinity on $\sigma$-compact manifolds (Theorem \ref{ind form}).

We remark that all the results presented here hold also for operators in the closure of the class of bounded \psos that are multiplication at infinity (Remark \ref{rmk ops cl}). On $\RR^n$, and more generally on manifolds that are asymptotically Euclidean, this closure contains the so-called isotropic 
 calculus, developed on $\RR^n$ by Shubin (\cite{shubin}) and in general by Melrose (\cite{melrose_scatt}). A thorough study of this enlarged class of operators, together with some consequences of the index formula, will appear in a forthcoming paper \cite{cc-nistor}. Moreover, the proof of cobordism invariance given here carries over, word by word, to the $G$-equivariant case, that is, to $G$-equivariant pseudodifferential operators, where $G$ is a (compact) Lie group. We also expect to obtain a similar result for families of pseudodifferential operators, as in \cite{a-sIV}.

\noindent\textit{Acknowledgments}. 
This work is part of a D Phil thesis, defended at the University of Oxford in 2003. I would like to thank my supervisors Prof Ulrike Tillmann and Prof Victor Nistor, from the Pennsylvania State University, for their encouragement and support, and  my thesis examinors for useful comments. I would also like to thank Klaas Landsman and Michael M\"uger for their helpful suggestions.

\section{The index}\label{s an-ind}

We consider a smooth, $\sigma$-compact\footnote{A manifold $M$ is said to be $\sigma$-compact if it is paracompact with a countable number of connected components; all our manifolds are assumed to be Hausdorff.}
 manifold $M$ without boundary,  smooth vector bundles $E, F$ over $M$, assumed to be trivial outside a compact set in $M$,\footnote{It suffices to assume that $E\cong F$ outside a compact set.} and define the class $\psimult(M;E,F)$ of \psos $$P: C^\infty_c(M;E)\to C^\infty(M;F)$$ that are multiplication at infinity. We show then that there is a well-defined index map $$\ind: K^0(TM)\to \ZZ$$ that computes the Fredholm index of elliptic operators in $\psimult(M;E,F)$.

\subsection{Pseudodifferential operators that are multiplication at infinity}\label{ss mpsos}

We first recall the main definitions of the theory of pseudodifferential operators, in order to fix notation (see \cite{horm-book, taylor, treves}).
Let $W\subset \RR^n$ be open; a smooth function $p:W\times\RR^n\into\CC$ is said belong to the 
class $S^m(W\times\RR^n)$ of \textit{symbols of order} $m$ if for any compact set $K\subset W$ and multi-indices $\alpha$, $\beta$ there exists $C_{K,\alpha,\beta} > 0$ such that
$$|\partial_x^\alpha \partial_\xi^\beta p(x,\xi)|\leq
C_{K,\alpha,\beta}(1+|\xi|)^{m-|\beta|} , $$
for all $x\in K$ and $\xi\in\RR^n$.
We always assume that the symbols are \textit{classical}, that is, that $p\in S^m(W\times\RR^n)$ has an asymptotic expansion $p\sim\sum {p_{m-k}}$
with $p_{m-k}\in S^{m-k}(W\times\RR^n)$
positively homogeneous of degree $m-k$ in $\xi$. 
%
%

Every element $p\in S^m(W\times\RR^n)$  induces a bounded operator
$p(x,D) : C^\infty_c (W) \rightarrow C^\infty (W) $
 given by 
\begin{equation}\label{pso RRn}
p(x,D)u(x) = (2\pi)^{-n}\int_{\RR^n}{p(x,\xi)\hat{u}(\xi)e^{ix\cdot \xi}
d\xi},
\end{equation}
where $\hat{u}$ denotes the Fourier transform of $u$.
 Such an operator $P=p(x,D)$ is called a \textit{\pso of order} (at most) $m$
 on $W$ with total symbol $p(x,\xi)$. 
The \textit{principal symbol} $\sigma_m(P)$ of $P=p(x,D)$ is  the class of $p(x,\xi)$ in $S^m(W\times\RR^n)\slash S^{m-1}(W\times\RR^n)$, that is, the leading term $p_m$ in the expansion of $p(x,\xi)$ as a classical symbol; it is a smooth function on $W\times \RR^n$,
positively homogeneous of degree $m$ in $\xi$ (i.e., for $ \| \xi \| > 1$ and $t\geq 1$, $p_{m}(x,t\xi) = t^{m}p_{m}(x,\xi)$).
%

The above definitions still make sense when $W\subset M$ is an open subset of a manifold $M$.  For operators on sections of vector bundles, locally,  we identify smooth sections with vector valued functions on the base space and an operator 
 $P: C^\infty_c(M;E)\to C^\infty(M;F)$ is said to be a \textit{pseudodifferential operator of order} (at most) 
$m$ on $M$, $P\in \Psi^m(M;E,F)$, if for any coordinate chart $W$ of $M$ trivializing $E$, $F$,  and for any $h\in C^\infty_c(W)$, $hPh: C_c^\infty(W)^N\to C_c^\infty(W)^N$ is a matrix of \psos of order $m$ on $W$. 

Regarding principal symbols, we have now for each $W$ a matrix of elements in $S^m(W\times\RR^n)\slash S^{m-1}(W\times\RR^n)$ and each of the entries of
this matrix transforms to define a smooth function on $T^*M $, the cotangent bundle of $M$. 
 The principal symbol is then globally defined as a smooth section of the bundle
$\Hom(E,F)$ over $T^*M$, 
 that is, it is a smooth bundle homomorphism
\begin{equation}
\sigma_m(P) : \pi^*E \to \pi^*E,
\end{equation}
 with $\pi: T^*M \to M$  the projection map, which is positively  homogeneous on the fibers of $T^*M$. The class of such symbols is identified with $$S^m(T^*M;E,F)\slash S^{m-1}(T^*M;E,F),$$ where $S^m(T^*M;E,F)$ is the class of smooth sections of $\Hom(E,F)$ over $T^*M$ whose pull-back over a coordinate chart $W$ gives a matrix of elements of $S^m(W\times\RR^n)$. 
 Since the equivalence class of a bundle morphism in $S^m(T^*M;E,F)\slash S^{m-1}(T^*M;E,F)$ does not depend on its behavior at $\xi=0$,   
  we can endow $M$ with a Riemannian metric and take its unit sphere bundle $S^*M$, to get that $S^m(T^*M;E,F)\slash S^{m-1}(T^*M;E,F)$ can be identified with
 a subclass of $C^\infty(S^*M;\Hom(E,F))$.


We can recover the principal symbol $\sigma_m(P)$ from a given \pso $P$ of order $m$: 
for $x_0\in M$, $\xi\in T^*M$, $\|\xi\|\geq 1$, and $u\in C^\infty_c(M;E)$,
\begin{equation}\label{princsymb}
\sigma_m(P)(x_0,\xi)u(x_0) =\lim_{t\rightarrow {+\infty}}t^{-m}\left({e^{- itf(x)} P(e^{itf(x)}u(x))}\right)(x_0),
\end{equation}
where $f\in C^\infty(M)$ is such that $df(x_0)=\xi$.

From now on, we assume that $M$ is Riemannian, endowed with a smooth measure, and $E$, $F$ have Hermitian metrics $\langle\,,\,\rangle_E$, $\langle\,,\,\rangle_F$. Denote by $(\,,\,)_E$ the usual inner product in $C^\infty_c(M;E)$ 
Recall that given a linear operator $P :C^\infty_c (M;E) \rightarrow C^\infty (M;F)$,
its \textit{
 adjoint} is defined as the linear operator $P^* :C^\infty_c (M;F) \rightarrow C^\infty (M;E) $ such that 
$$( Pu,v)_F = ( u,P^*v)_E, $$
for all $u\in C^\infty_c(M;E)$, $v\in C^\infty_c(M;F)$; if $P$ is a \pso of order $m$ then $P^*$ also is and
 $\sigma_m(P^*)=\sigma_m(P)^*.$
%

To define the class of operators we will be working with, we first give the following definition: a linear operator $P : C^\infty_c(M;E)\to C^\infty(M;F)$ is said to be \textit{multiplication at infinity} if there is a compact $K\subset M$  and a smooth bundle homomorphism $p : E\rightarrow F$ such that, defining $C^\infty_c(M-K;E):=\{u\in C^\infty_c(M;E): \supp u \subset M-K\}$, 
\begin{equation}\label{ops mult infty}
Pu = pu,\;\text{for}\; u\in C^\infty_c(M-K;E).
\end{equation}
\begin{defn}\label{mpsos}
A  \pso  $P :C^\infty_c (M;E) \rightarrow C^\infty (M;F) $ of order $m$ is said to be in the class $\Psi_{mult}^m(M;E,F)$ if
$P$ and $P^*$ are multiplication operators at infinity.
\end{defn}
 
 Clearly, if $M$ is compact, $\Psi_{mult}^m(M;E,F)$ coincides with the usual class of \psos of order $m$. 

If $P : C^\infty_c(M;E)\to C^\infty(M;F)$ is such that there is a morphism $p$ such that $Pu=pu$ outside some compact $K$, 
then (\ref{princsymb}) yields that, for $x\notin K$, $\|\xi\|\geq 1$, $\sigma_m(P)(x,\xi)=0$ for $m>0$ and $\sigma_m(P)(x,\xi)=p_x$ for $m\leq 0$ (it is clear that for $m<0$, $p=0$). In particular, if $P$ is multiplication at infinity, then the principal symbol is constant on the fibers of $T^ *M$ outside a compact set in $M$. 
 The (partial) converse can be given as follows: if $P$ is a \pso on $M$ and if $\sigma_m(x,\xi)$ does not depend on $\xi$, for $x\notin K\subset M$ compact, then $P$ is multiplication at infinity modulo operators of lower order.
(In particular, if $P$ is a \pso that is multiplication at infinity, then $P^*$ is multiplication at infinity modulo \psos of lower order.)

The following proposition gives an equivalent definition of  $\Psi_{mult}^m(M;E,F)$; let $C^\infty_K(M;E)$ denote the class of smooth, compactly supported sections of $E$ whose support is contained in some fixed compact $K\subset M$.

\begin{prop}\label{P=P_1+p}
Let $P:C^\infty_c (M;E) \rightarrow C^\infty (M;F) $ be a \pso of order $m$. Then $P\in\Psi_{mult}^m(M;E,F)$ if and only if 
$$P=P_1+p,$$
for a smooth bundle morphism $p:E\to F$, with $p=0$ for $m<0$, and an order $m$ \pso $P_1 :C^\infty_c (M;E) \rightarrow C^\infty(M;F)$ such that $\Img P_1\subset C^\infty_K(M;F)$  and $\Img P_1^*\subset C^\infty_K(M;E)$  for some compact $K\subset M$.
\end{prop}
\pf
Let $P\in\Psi_{mult}^m(M;E)$, $p:E\to E$ and $K\subset M$ compact be such that $Pu=pu$ outside $K$ (note that we have then $p=0$, for $m<0$). Define $P_1 = P -p$. Then  $P_1$ is a \pso of order $m$, since $P$ is and $p$ has order zero\footnote{We denote the morphism $p$ and the operator of multiplication it induces by the same letter $p$.}. For all $v\in C^\infty_c (M;F)$ such that $\supp v\subset M-K$, we have $P_1^*v = P^*v -p^*v = 0$. 
But then one has also that, for all $u\in C^\infty_c(M;E)$,
$$(P_1u,v)_E = ( u,P_1^*v)_E = 0$$
 and hence  $P_1u(x)=0$, for all $x\notin K$. Using the same reasoning for the adjoint $P_1^*$, we have then that $\supp(Pu)\subset K$, for any $u\in C_c^\infty(M;E)$ and  $\supp(P^*v)\subset K$, for any $v\in C_c^\infty(M;E)$. 

Now let $P$ be such that there exist $P_1$, $K\subset M$ compact and $p:E\to F$ in the conditions above, with $p=0$, for $m<0$. For any  $u\in C^\infty_c(M;E)$ such that $ \supp u\subset M-K$, and $v\in C^\infty_c(M;E)$, one has
$$( P_1u,v)_E = ( u,P_1^*v)_E = 0,$$
since $\supp u\cap \supp(P_1^*v)=\emptyset$. We conclude that $P_1u =0$, that is, that $Pu=pu$ and $P$ is multiplication at infinity.  
Again, the same reasoning applies to the adjoint $P^*$ and therefore $P\in\Psi_{mult}^m(M;E,F)$. 
\qed

Note that a \pso $P_1$ satisfies the condition in the previous proposition if and only if its distribution kernel has compact support contained in $K\times K$. 

The above result allows the generalization of the main results of the theory of \psos on compact manifolds to $\Psi^m_{mult}(M;E,F)$. We give here a only a summary of the relevant results (a detailed account will be given in \cite{cc-nistor}). 

First note that if $P$ is multiplication at infinity, then for $u\in C^\infty_c(M;E)$,  we have $\supp(Pu)\cap (M- K)=\supp u\cap(M-K)$, and hence, 
\begin{equation}\label{properlysupp}
\supp (Pu)\subset K\cup \supp u.
\end{equation}
 We conclude that $P$ and $P^*$ map compactly supported sections to compactly supported sections, that is, $P$ is properly supported.  
 Denote by  $S^m_{mult}(T^*M;E,F)$ the set of equivalence classes of $p(x,\xi)$ in $S^m(T^*M;E,F)\slash S^{m-1}(T^*M;E,F)$ such that $p(x,\xi)$ is independent of $\xi$, for $x\notin K$, $K\subset M$ compact. The principal symbol map then gives rise to a map 
\begin{equation}\label{symb map mpsos}
\sigma_m: \Psi_{mult}^m(M;E,F)\to S^m_{mult}(T^*M;E,F).
\end{equation}
We have that, if $P\in\Psi_{mult}^m(M;E,F)$, $Q\in\Psi_{mult}^l(M;F,G)$, then $P^*\in\Psi_{mult}^m(M;F,E)$, $QP\in P\in\Psi_{mult}^{m+l}(M;E,G)$ with
$$\sigma_m(P^*)=\sigma_m(P)^*,\;\;\sigma_{m+l}(QP)=\sigma_l(Q)\,\sigma_m(P).$$
In particular, considering $E=F$, we have that  $\Psi_{mult}^0(M,E)$ is a $*$-algebra and $\sigma_0$ is a $*$-homomorphism.
Moreover, one also has, as in the compact case, that the principal symbol map $\sigma_m$ is surjective and has kernel $\Psi^{m-1}_{mult}(M;E,F)$.

From now on, we will be mostly concerned with operators of order (at most) zero. When $M$ is a \textit{compact} manifold, \psos of order zero on $M$ (acting on sections of two vector bundles $E$ and $F$) can  be extended as bounded operators  $L^2(M;E)\to L^2(M;F)$, where $L^2(M;E)$ the completion of $C_c^\infty(M;E)$ with respect to the norm
$$\| u\|_{L^2(M,E)} := \left({\int_M\langle u(x),u(x)\rangle_E}\right)^{\frac{1}{2}}.$$
This result can be generalized to \psos with compactly supported kernel, like the ones appearing in Proposition \ref{P=P_1+p}. The following result is therefore not surprising.

\begin{prop}\label{psos bdd}
Let  $P\in\psimult(M;E,F)$. Then $P$ can be extended as a bounded operator $P: L^2(M;E)\to L^2(M;F)$ if, and only if, its principal symbol $\sigma_0(P)(x,\xi)\in S^0_{mult}(T^*M;E,F)$ is bounded in the $\sup$ norm, as a section of $\Hom(E,F)$ over $T^*M$. 
\qed
\end{prop}
%
One can also show, from the compact case, that a \pso with negative order that is given by a compactly supported kernel is necessarily compact. Proposition \ref{P=P_1+p} now yields:

\begin{prop}\label{psos compact}
Let  $P\in\Psi_{mult}^m(M;E,F)$, with $m<0$. 
 Then $P: L^2(M;E)\to L^2(M;F)$ is a compact operator. 
\end{prop} 
\pf 
From Proposition \ref{P=P_1+p}, there is a  compact set $K\subset M$ such that $\Img(P)\subset C^\infty_K(M;E)$ (and the same holds for $P^*$). The proof now goes as in the compact case.
\qed

\subsection{The index map}\label{ss ind map}

We will restrict the study of the Fredholm index to bounded operators (this is always possible by replacing an unbounded $P$ by $P(1+P^*P)^{-1\slash 2}$). From now on, we consider the class
\begin{equation}\label{def psimult +bdd}
\Psi_{mult}^{b}(M;E,F):=\{P\in\psimult(M;E,F): \sigma_0(P)\;\text{is bounded}\}.
\end{equation}

\begin{defn}\label{ell op}
An operator $P\in\Psi_{mult}^b(M;E,F)$ is said to be \textup{bounded elliptic} if its principal symbol can be represented by $p(x,\xi)\in S^0_{mult}(M;E,F)$ such that $p(x,\xi)$ is invertible in $C_b(S^*M;E,F)$.
\end{defn}
The above condition means that $p(x,\xi):E_x\to F_x$ is an isomorphism for $\xi\neq 0$ and  $p^{-1}(x,\xi)$ is bounded. 
In this case, $p^{-1}(x,\xi)$ is always positively homogeneous of degree zero and constant on the fibers outside some compact subset of $M$, that is, $p^{-1}(x,\xi)$ defines a class in $S^0_{mult}(T^*M;F,E)$.

\begin{prop}\label{fred}
Let $P\in \Psi_{mult}^b(M;E,F)$ be bounded elliptic. Then $P$ is a Fredholm operator.
\end{prop}
\pf
Since $P$ is bounded, it suffices to show that there exists $Q\in{\Psi}_{mult}^b(M;F,E)$ such that $PQ - I_F\in\K(M;F)$ and $QP - I_E\in\K(M;E)$, where $\K(M;E)$ denotes the class of compact operators $L^2(M;E)\to L^2(M;E)$.  
Let $p(x,\xi) := \sigma_0(P)(x,\xi)$ be invertible for $\xi\neq 0$ and $q(x,\xi):= p^{-1}(x,\xi)$. Since the principal symbol map is surjective, there exists  $Q\in \Psi^0_{mult}(M;F,E)$ such that $\sigma_0(Q)(x,\xi)=q(x,\xi)$, and $Q$ is bounded since $\sigma_0(Q)$ is. In this case, we have 
$$\sigma_0(PQ - I_F) = \sigma_0(P)\sigma_0(Q) - I = 0,$$
and similarly for $QP-I_E$, which yields  that $PQ - I_F\in \Psi^{-1}_{mult}(M;F)$ and $QP-I_E\in \Psi^{-1}_{mult}(M;E)$. From Proposition \ref{psos compact}, the result follows.
\qed

One can therefore consider, for every bounded elliptic operator $P\in \Psi_{mult}(M;E,F)$, the \textit{Fredholm index}
\begin{equation}\label{fred ind}
\ind (P) = \dim \ker P - \dim \coker P\in\,\ZZ.
\end{equation}
Since $P$ is bounded (and has closed range, since it is Fredholm), we have $\dim \coker P = \dim \ker P^*$.
As in the compact case, the Fredholm index of a bounded elliptic operator $P\in \Psi_{mult}^b(M;E,F)$ only depends on the homotopy class of $\sigma_0(P)$, and, as we will see now, it can be expressed as a $K$-theory map.

For the details on the following contructions, see \cite{atiyah-K}.  
For a locally space $X$, let $K^0(X)$ denote the compactly supported $K$-theory
of $X$. Elements in $K^0(X)$ are given by equivalence classes of triples $(E,F,\alpha)$, where $E,F$ are vector bundles over $X$ and $\alpha:E\to F$ is invertible outside a compact subset of $X$; two triples $(E,F,\alpha)$ and $(E^\prime,F^\prime,\alpha^\prime)$ are said to be equivalent if there exist vector bundles $G$, $H$ such that $E\oplus G\cong E^\prime\oplus H$, $F\oplus G\cong F^\prime \oplus H$, and these isomorphisms are compatible with $\alpha, \alpha^\prime$. We denote by $[E,F,\alpha]$ the equivalence class of $(E,F,\alpha)$. 

Let $P\in\Psi_{mult}^b(M;E,F)$ be a bounded elliptic operator, with principal symbol 
$$\sigma_0(P)(x,\xi)= p(x,\xi)\in S^0_{mult}(T^*M;E,F)\subset C^\infty(T^*M;\Hom(E,F));$$ 
then $p(x,\xi)$ is an isomorphism for $\|\xi\| \neq 0$. On the other hand, since $P$ is multiplication at infinity, there is a compact $K\subset M$ such that $p(x,\xi)$ does not depend on $\xi$, for $x\notin K$. Hence, for $x\notin K$, we have that $p(x,\xi)$ is an isomorphism for all $\xi$. 
 Regarding $p(x,\xi)$ as a bundle map $\pi^*E\to \pi^*F$, $\pi:T^*M\to M$, we have that the triple $(\pi^*E,\pi^*F, p)$ defines a $K$-theory class  
\begin{equation}\label{sigma_0 in K(TM)}
[\sigma_0(P)]:= [\pi^*E,\pi^*F, p]\in K^0(T^*M),
\end{equation} 
called the \textit{symbol class} of $P$. 

Using the equivalence relation given above, it is not hard to show that, as in the compact case, $\ind (P)$ only depends on the symbol class $[\sigma_0(P)]$, and we would like to have the index as a $K$-theory map. To see that that is indeed the case, we make use of the following result (the proof is just a little adaptation of the argument given in \cite{a-sI}, page 492).\footnote{Note that, since $M$ is a manifold, it has finite topological dimension, and every vector bundle over $M$ can be complemented to a trivial bundle; see \cite{hus}, page 31.}

\begin{lem}\label{lemma K-thry = symb}
Let $V$ be a vector bundle over a manifold $M$. Every element in $K^0(V)$ can be written as a class $[\pi^*E,\pi^*F,\alpha]$, where $E, F$ are vector bundles over $M$ that are trivial outside a compact, and $\alpha:\pi^*E\to\pi^*F$ is 
 positively homogeneous of degree zero and constant on the fibers outside a compact $K\subset M$. Moreover, we can assume that $E_{|M-K}=F_{|M-K}=\theta^n$ and that $\alpha_{|V-\pi^{-1}(K)}=1$ (where $\theta^n$ denotes the trivial bundle). 
\qed
\end{lem}
Note that the morphism $\alpha$ given in the above lemma is always bounded, since it coincides with the identity over $V-\pi^{-1}(K)$ and it is positively homogeneous of degree zero,<and the same argument shows that its inverse, defined outside a compact in $V$, is also bounded. Using the fact that the principal symbol map is surjective, one finally gets:
\begin{prop}\label{K-thry = symb}
Every element in $K^0(T^*M)$ is the symbol class of a bounded elliptic operator in $\Psi_{mult}^b(M;E,F)$.
\qed
\end{prop}
Using the metric  to identify canonically $T^*M$ with $TM$, we have then that there is a well-defined map 
\begin{equation}\label{an-ind}
\ind : K^0(TM) \to \ZZ
\end{equation}
such that, for any bounded elliptic operator $P\in\Psi_{mult}^b(M;E,F)$,
\begin{equation}\label{anind sigma=ind P}
\ind[\sigma_0(P)] = \ind(P)= \dim\left(\ker(P)\right)-\dim\left(\ker(P^*)\right).
\end{equation}
One important property of this map, which acquires new features in the non-compact case, is its invariance with respect to open embeddings, the so-called \textit{excision property}, given by Atiyah and Singer in \cite{a-sI} for compact manifolds. The proof holds with small modifications also for non-compact $M$, where now we compute the indices of operators that are multiplication at infinity.

Recall that for a locally compact space $X$ and open $U\subset X$, there is a natural
 map $j:X^+\to X^+/(X^+-U)=U^+$ which induces a $K$-theory map 
\begin{equation}\label{ext}
j_*: K^0(U)\to K^0(X).
\end{equation}
In particular, if $M$ is a manifold and $U\subset M$ is open, then $TM_{|U}=TU$ is also open in $TM$ and we have 
\begin{equation}\label{ext TX}
h_{M} : K^0(TU) \to K^0(TM).
\end{equation} 

\begin{prop}\label{an-ind vs ext}
Let $M$ be a manifold, $U\subset M$ be open and $$h_M: K^0(TU)\to K^0(TM)$$ be the induced map as in (\ref{ext TX}). Then $\ind \circ h_M =\ind$.
\end{prop}
\pf
Given $a\in K^0(TU)$, write it first as the principal symbol of a bounded elliptic operator $P\in\Psi_{mult}^b(U;E,F)$, $a=[\pi^*E,\pi^*F,p]$, with $E$, $F$ and $p(x,\xi):=\sigma_0(P)(x,\xi)$ satisfying the extra conditions given in Lemma \ref{lemma K-thry = symb}.
 Then $E$, $F$ can be extended as vector bundles $\tilde{E},\tilde{F}$ over $M$, which are trivial in $M\backslash U$, and, similarly, $p(x,\xi)$ can be extended as a map $\tilde{p}: \pi_M^*\tilde{E}\to\pi_M^*\tilde{F}$, which is the identity on $TM\backslash TU$.  By carefully looking at the map $h_M$, one checks that
$$h_M(a)=h_M[\pi^*E,\pi^*F,p]= [\pi_M^*\tilde{E},\pi_M^*\tilde{F},\tilde{p}].$$
Now write $P=P_1+f$, as in Proposition \ref{P=P_1+p}. 
Since $P$ is multiplication outside a compact $K\subset U\subset M$, one can easily define an extension of $P_1$ to a \pso $\tilde{P_1}: C_c^\infty(M;\tilde{E})\to C^\infty(M;\tilde{F})$ such that the kernel of $\tilde{P_1}$ has support in $K\times K\subset M\times M$; on the other hand, since $f_x=p_x$ for $x\notin K$, we can also extend $f$ to $\tilde{f}:\tilde{E}\to \tilde{F}$ using $\tilde{p}$. 
We can then define an operator $\tilde{P}\in\Psi_{mult}^b(M;\tilde{E},\tilde{F})$ by
$$\tilde{P}:= \tilde{P}_1 + \tilde{p}$$
and clearly $[\sigma_0(\tilde{P})]=h_M(a)$. It is easy to check that, by construction, $\ker(P)=\ker(\tilde{P})$, $\coker(P)=\coker(\tilde{P})$. Hence the Fredholm indices of $P$ and $\tilde{P}$ coincide and the result follows.
\qed

We finish the section with the following remark.
\begin{rmk}[The closure of $\Psi_{mult}^b(M;E,F)$]\label{rmk ops cl}
\upshape
Let
$$\overline{\Psi^b}_{mult}(M;E,F):=\overline{\Psi_{mult}^b(M;E,F)},$$
where the closure is, of course, taken with respect to the operator norm.
 It follows from (\ref{princsymb}) that the principal symbol map $$\sigma_0 :\psimult(M;E,F) \to S^0_{mult}(T^*M;E,F)$$
 is continuous when restricted to $\Psi_{mult}^b(M;E,F)$, where we endow the subclass of $S^0_{mult}(T^*M;E,F)$ of bounded symbols with the $\sup$-norm.
 The closure of this class of bounded symbols is seen to be given by the class $C_{a}(S^*M;E,F)\subset C_b(S^*M;\Hom(E,F))$ of continuous sections of $S^*M$ that are 'asymptotically constant on the fibers'. We therefore get an extended symbol map
$${\sigma_e}: \overline{\Psi^b}_{mult}(M;E,F) \to  C_{a}(S^*M;E,F),$$
which fits into the  exact sequence of $C^*$-algebras
\begin{equation*}
\begin{CD}
0 @>>> \K(M;E) @>>> \overline{\Psi^b}_ {mult}(M;E) @>{\sigma_e}>>C_{a}(S^*M;E) @>>> 0,
\end{CD}
\end{equation*}
where $\K(M;E)$ denotes the class of compact operators $L^2(M;E)\to L^2(M;E)$. An operator $P\in  \overline{\Psi^b}_{mult}(M;E,F)$ whose extended symbol $\sigma_e(P)$ is invertible in $C_b(S^*M;E,F)$ is a Fredholm operator, and one can define a symbol class $[\sigma_e(P)]\in K^0(TM)$ such that the index map (\ref{an-ind}) also computes the index in $\overline{\Psi^b}_ {mult}(M;E,F)$. (See \cite{cc-nistor} for the proofs of these results.)

All the results given here for $\Psi_{mult}^b(M;E,F)$ will hold more generally in $\overline{\Psi^b}_{mult}(M;E,F)$, replacing $\sigma_0$ by $\sigma_e$. Note that this class contains very familiar classes of operators. For instance, when $M=\RR^n$ and $E=F=\CC$, we have that
$$\overline{\Psi^b}_{mult}(\RR^n)\supset\Psi_{iso}(\RR^n)$$
where $\Psi_{iso}(\RR^n)$ denotes the so-called isotropic 
 calculus on $\RR^n$, developed extensively, among others, by Shubin (see \cite{melrose_scatt, shubin}).
\end{rmk}

\section{Cobordism of symbols}\label{s cobd inv}

We will define now a relation of cobordism for pairs $(M,\sigma)$, where $M$ is, as usual, a $\sigma$-compact manifold (without boundary) and $\sigma\in K^0(TM)$. We start with defining the map of restriction of symbols. First recall that, for any locally compact space $X$, Bott periodicity yields an isomorphism
\begin{equation}\label{bott iso}
\beta_X: K^0(X)\to K^0(X\times\RR^2),\; a\mapsto a\cup \beta_{\RR^2}
\end{equation}
called the \textit{Bott isomorphism}, where $\beta_{\RR^2}:=\beta_{\RR^2}(1)\in K^0(\RR^2)$ is the so-called Bott class.

Now let $X$ be a manifold with boundary and assume that $M$ is a boundary component of $X$, that is, $M$ is a connected submanifold of $\partial X$.
Using Bott periodicity, we can identify $\sigma\in
K^0(TM)$ with $\beta_{TM}(\sigma) \in K^0(TM \times \RR^2)$, where $\beta_{TM}: K^0(TM)\to K^0(TM\times\RR^2)$ denotes the Bott
isomorphism. 
On the other hand, since $ M \subset \partial X$, the normal bundle to $M$ in $X$ is the trivial one-dimensional vector bundle,
 and we have
\begin{equation}\label{TX|bdry inwd}
TX_{| M} \cong TM \times \RR,
\end{equation}
where we always assume that the isomorphism is given by taking the \textit{inward} normal vector at each point of $M$.
Now, since $TX_{| M}$ is closed in $TX$, we have a restriction map  $K^1(TX) \rightarrow K^1(TX_{|M})$, induced by the inclusion $TX_{| M} \times \RR \to TX \times \RR$ and, using the isomorphism (\ref{TX|bdry inwd}), we now get a map 
\begin{equation}\label{r_M}
r_M : K^1(TX) \rightarrow K^1(TM \times \RR):=K^0(TM\times\RR^2).\footnote{Note that if we consider in (\ref{TX|bdry inwd}) the isomorphism given by the outward normal vector, we  get $-r_M$.}
\end{equation}
 Define the \textit{symbol restriction} map as
\begin{equation}\label{u_M}
u_M = \beta_{TM}^{-1} \circ r_M : K^1(TX) \rightarrow K^0(TM).
\end{equation}
We can now make clear what we mean with \textit{cobordism of symbols}.

\begin{defn}\label{codsm symb def}
Let $M_1$, $M_2$ be manifolds, and $\sigma_i\in K^0(TM_i)$, $i=1,2$. We say that the pairs $(M_1,\sigma_1)$ and $(M_2,\sigma_2)$ are \textup{symbol-cobordant} if there exists a manifold $X$, and $\omega \in K^1(TX)$ such that

(i) $\partial X= M_1 \sqcup M_2$

(ii) $u_{M_1}(\omega)=\sigma_1$, and  $u_{M_2}(\omega)=-\sigma_2$.

\noindent We write  $(M_1,\sigma_1)\sim (M_2,\sigma_2)$ and the pair $(X,\omega)$ is called a \textup{cobordism of symbols}.
\end{defn}
%

Note that the map $r_M$ is given by the restriction map in the long exact sequence 
\begin{equation*}
\begin{CD}
... @>>> K^1(TX)@>{r_M}>>K^1(TX_{|M})  @>{\partial_i}>>K^0(TX,TX_{|M}) @>>> ...,
\end{CD}
\end{equation*}
and condition $(ii)$ is equivalent to $\partial_i(\beta_{TM_i}(\sigma_i))=0$, where $\partial_i$ denotes the connecting map.

We define an operation of direct sum of pairs, for manifolds with the same dimension, by
\begin{equation}\label{dir sum symb}
(M_1,\sigma_1)\oplus (M_2,\sigma_2):= (M_1 \sqcup M_2, \sigma_1\oplus \sigma_2),
\end{equation}
where $\sqcup$ denotes disjoint union (and we are using the fact that $K^0(T(M_1 \sqcup M_2))=K^0(TM_1)\oplus K^0(TM_2)$).
 It is easy to check that, for each $n\in\NN$, cobordism of symbols defines an equivalence relation on $$\{(M,\sigma): M\;\text{is an open manifold,}\;\dim M =n, \sigma\in K^0(TM)\}$$
 and the set of equivalence classes is an abelian group with the direct sum defined by (\ref{dir sum symb}).

\begin{rmk}\label{rmk cbd k-thry}
\upshape
The notion of cobordism of symbols is closely related to the notion of cobordism in $K$-theory, as developed in \cite{palais} for \textit{compact} manifolds: two classes $a_i\in K^0(M_i)$, $i=1,2$ are said to be cobordant if there is a compact manifold $X$ such  that $\partial X=M_1\sqcup M_2$ and $b\in K^0(TX)$ such that $b_{|M_i}=a_i$. This group arose in relation to cobordism for (symbol classes of) twisted signature operators in \cite{a-s-cobd}. 
%
%
%
\end{rmk}

In the following proposition we check that cobordism of symbols is stable under the cup product with an element of $K^0(M)$, more precisely, with a class in the cobordism group mentioned in Remark \ref{rmk cbd k-thry},  as long as there exists a common cobordism manifold. (If $a\in K^0(X)$, $b\in K^0(Y)$, for locally compact space $X, Y$, then $a\cup b$ is defined as $a\cup b:=\pi_X^*a\otimes \pi_Y^*b\in K^0(X\times Y)$, with $\pi_X:X\times Y\to X$, $\pi_Y:X\times Y\to Y$ the projection maps.)

\begin{prop}\label{twist symb}
Let $M_i$ be a manifold and $\sigma_i\in K^0(TM_i)$, $i=1,2$. If  $(M_1,\sigma_1)\sim (M_2,\sigma_2)$, with $(X,\omega)$ giving a cobordism of symbols, and if $a_i\in K^0(M_i)$, $i=1,2$, are such that there exists $b\in K^0(X)$ with $b_{|M_i}=a_i$ then  $(M_1,\sigma_1\cup a_1)\sim (M_2,\sigma_2\cup a_2)$.
\end{prop}
\pf 
It suffices to show that if  $(M,\sigma)\sim 0$, with $(X,\omega)$ a cobordism, and if $a\in K^0(M)$ is such that there is $b\in K^0(X)$ with $b_{|M}=a$, then $(M,\sigma\cup a)\sim 0$. 
 Let $\beta_{TM} : K^0(TM) \to K^0(TM\times\RR^2)$ be the Bott isomorphism and $r_M : K^1(TX) \to K^1(TM\times \RR)$ be the restriction map as in (\ref{r_M}).
 Then $r_M(\omega) = \beta_{TM}(\sigma)$ and
\begin{eqnarray*}
\beta_{TM}(\sigma\cup a) &  = & \beta_{\RR^2}\cup(\sigma\otimes\pi_{TM}^*(a))\\
& = & \beta_{TM}(\sigma)\otimes \pi_1^*\circ\pi^*_{TM}(a) \\
& = & r_{M}(\omega)\otimes(\pi_{TM}\circ\pi_1)^*(a), 
\end{eqnarray*}
where $\pi_{TM} : TM \to M$, $\pi_1: TM\times \RR^2 \to TM$, $\pi_2 : TM \times \RR^2 \to \RR^2$ are the projection maps, and $\beta_{\RR^2} \in K^0(\RR^2)$ is the Bott class. 
Writing $a=b_{|M}= i^*(b)$, with $i:M\to X$ the inclusion, and $r_M = k^*$, where again $k : TX_{|M}\times\RR \to TX\times\RR$ is the inclusion, and noting that $TX_{|M}\times\RR\cong TM\times\RR^2$, we have $\pi_{TX}\circ k = i\circ \pi_{TM}\circ\pi_1$, with $\pi_{TX} :TX\times \RR \to X$ the projection map. Hence, 
$$\beta_{TM}(\sigma\cup a) = r_M(\omega) \otimes r_M(\pi^*_{TX}(b))=r_M(\omega\cup b).$$
We conclude that $(M,\sigma\cup a)\sim 0$ via the cobordism $(X,\tilde{\omega})$, with $\tilde{\omega}=\omega\cup b$. 

\qed

Recall that the \textit{signature operator} on a manifold $M$, oriented and even-dimensional, is given by 
\begin{equation}
D_M=d+d^*: C^\infty(M; \Lambda^+(T^*M)) \to C^\infty(M; \Lambda^-(T^*M)), 
\end{equation}
where $d$ is the exterior derivative, and the $\ZZ_2$-grading is given by the chirality operator $\Gamma:= i^{n\slash 2}e_1...e_n$, with $e_1,...,e_n$ an oriented basis for $T_xM$, $x\in M$ (see \cite{lawson-mich}).
 Its symbol is given by the triple
\begin{equation}\label{symb sign}
  \sigma(D_M)=(\pi^*\Lambda^+(T^*M), \pi^*\Lambda^-(T^*M), c),
\end{equation}
where $\pi: T^*X \to X$ and $c$ is given by Clifford multiplication, 
 that is, $c{(\xi)}\alpha = \xi\wedge\alpha \;-\; i_\xi(\alpha)$, for $\alpha\in\Lambda^+(T_x^*X)$, $\xi \in T^*_xX$. If $M$ is compact, then $\sigma(D_M)$ defines a class $[\sigma(D_M)] \in K^0(T^*M)$.\footnote{In general, one has $[\sigma(D_M)] \in KK(C_0(M),C_0(T^*M))$.} 
Given a vector bundle $V$ over $M$, one can similarly define the so-called \textit{twisted signature operator}
\begin{equation}
D_W: C^\infty(M; \Lambda^+(T^*M;W)) \to C^\infty(X; \Lambda^-(T^*M;W)), 
\end{equation}
and, when $M$ is compact, its symbol is given by $\sigma_V=\sigma_M\cup [W]$, where $[W]\in K^0(M)$ denotes the equivalence class defined by $W$. 
 We will make use of the following lemma, whose proof is straightforward (it is basically a consequence of the multiplicativity property of Thom classes).

\begin{lem}\label{lem symb sign}
(i) Let $\beta_{\RR^2}\in K^0(\RR^2)$ be the Bott class and let  $\sigma_{\RR^2}$ denote the symbol  of the signature operator on $\RR^2$. Then
\begin{equation*}\label{sigmarr2=2beta} 
\sigma_{\RR^2} = \pi^*\beta_{\RR^2} \oplus\pi^*\beta_{\RR^2},
\end{equation*}
where $\pi: T\RR^2\to\RR^2$ is the projection map. 

(ii) Let $D_M, D_{M^\prime}$ be the signature operators on even-dimensional, oriented manifolds $M, M^\prime$. If $M$ is compact, then  
\begin{equation*}\label{sigma XtimesY=sigmaxtimessigmaY}
\sigma(D_{M\times M^\prime})=\sigma(D_M)\cup\sigma(D_{M^\prime})\in K^0(T(M\times M^\prime)).
\end{equation*}
\qed
\end{lem}
 
\begin{prop}\label{sign sim 0}
Let $M$ be a compact, even-dimensional, oriented manifold, and $D_M$ the signature operator on $M$, with principal symbol $\sigma_M\in K^0(TM)$. If there is a  compact manifold $X$ such that $\partial X = M$, then $(M,\sigma_M\oplus\sigma_M)\sim 0$.
\end{prop}
\pf
Consider the signature operator on $X\times\RR$, with symbol $\sigma_X$. Restricting to $X\times\{0\}$, we can define $$\omega:= (\sigma_X)_{|TX\times \RR}\in K^0(TX\times\RR)=K^1(TX).$$  
Now consider the signature operator on $M\times\RR^2$, with symbol $\sigma_{M\times\RR^2}$ and restrict to $M\times\{0\}$. One can check directly that
$$r_M(\omega)=(\sigma_{M\times\RR^2})_{|M\times\{0\}}$$
where $r_M: K^1(TX)\to K^1(TX_{|M})=K^0(TM\times\RR^2)$ is the restriction map.
From Lemma \ref{lem symb sign}, we have 
$$\sigma_{M\times\RR^2}=\sigma_M\cup\sigma_{\RR^2}=(\sigma_M\cup\pi^*\beta_{\RR^2})\oplus(\sigma_M\cup\pi^*\beta_{\RR^2}).$$
Therefore, 
$$r_M(\omega)=(\sigma_M\cup\beta_{\RR^2})\oplus(\sigma_M\cup\beta_{\RR^2})=\beta_{TM}(\sigma_M\oplus\sigma_M),$$
and then $(M,\sigma_M\oplus\sigma_M)\sim 0$.
\qed

We now finally state our main result on cobordism invariance. 

\begin{thm}\label{cobdsm inv}
Let $M$ be a  $\sigma$-compact manifold and $P\in\Psi_{mult}^b(M;E,F)$ a bounded elliptic operator. Let $\sigma:=[\sigma_0(P)] \in K^0(TM)$ denote the principal symbol class of $P$. If $(M,\sigma)\sim 0$, then $\ind(P)=0$.
\end{thm}
We postpone the proof of this theorem until Section \ref{s pf cbd inv} and now give some corollaries. As before, $\beta_{TM}: K^0(TM)\to K^0(TM\times\RR^2)$ and $r_M: K^1(TX)\to K^1(TM\times\RR)$ denote the Bott isomorphism and the restriction map, respectively, and  $\partial : K^1(TX_{|M})\to K^0(TX,TX_{|M})$ denotes the connecting map.

\begin{cor}\label{cobd inv cmpt}
Let $M$ be a compact manifold,  and let $P\in\Psi^0(M;E,F)$ be an elliptic \pso with principal symbol $\sigma=[\sigma_0(P)]\in K^0(TM)$. Assume that there is a compact manifold $X$ such that $\partial X=M$; if there exists $\omega \in K^1(TX)$ such that $r_M(\omega)=\beta_{TM}(\sigma)$, or equivalently, if $\partial(\beta_{TM}(\sigma))=0$, then $\ind(P)=0$.
\end{cor}
\pf
Straighforward from Theorem \ref{cobdsm inv} and the definition of cobordism of symbols.
\qed

The following corollary gives the cobordism invariance of the index of signature operators on compact manifolds, as in \cite{a-s-cobd, palais}.

\begin{cor}\label{cobd inv sign ops}
Let $M$ be a compact, even-dimensional, oriented manifold, $V$ a vector bundle over $M$, and $D_V$ the twisted signature operator on $M$. If there is a compact manifold $X$ with $\partial X = M$ and a vector bundle $W$ over $X$ such that $W_{|M}=V$, then $\ind(D_V)=0$.
\end{cor}
\pf From Theorem \ref{cobdsm inv} and Propositions  \ref{twist symb}, \ref{sign sim 0}, we have that $$\ind(D_V\oplus D_V)=0.$$ The additivity of the Fredholm index gives the result.
\qed

To conclude, we check that Theorem \ref{cobdsm inv} yields that in fact the index map is well defined on cobordism classes.

\begin{cor}\label{cobds inv sum}
Let $P\in{\Psi}^b_{mult}(M;E,F)$, $P^\prime\in{\Psi}^b_{mult}(M^\prime;E^\prime,F^\prime)$ be bounded elliptic, and denote by $\sigma:=[\sigma_0(P)]\in K^0(TM)$, $\sigma^\prime:=[\sigma_0(P^\prime)]\in K^0(TM^\prime)$ the principal symbol classes. If $(M,\sigma)\sim (M^\prime,\sigma^\prime)$, then $\ind(P)=\ind(P^\prime)$.
\end{cor}
\pf
Let $P^*\in{\Psi}_{mult}(M;F,E)$ be the adjoint of $P$; then $P^*$ is Fredholm, with $\ind(P^*)= -\ind(P)$. We have also that $\sigma_0(P^*)=(\sigma_0(P))^*$ and therefore $[\sigma_0(P^*)]=-[\sigma_0(P)]=-\sigma\in K^0(TM)$.
 Now, since $(M,\sigma)\sim (M^\prime,\sigma^\prime)$, we have $(M\sqcup M^\prime,-\sigma \oplus \sigma^\prime)\sim 0$. Consider  $P^*\oplus P^\prime: C^\infty_c(M\sqcup M^\prime;F\sqcup E^\prime)\to C^\infty_c(M\sqcup M^\prime;E\sqcup F^\prime)$; then  $-\sigma \oplus \sigma^\prime = \sigma_0(P^*\oplus P^\prime)$ 
and,  from Theorem \ref{cobdsm inv}, we conclude that  $$\ind(P^*\oplus P^\prime)=\ind(P^*)+\ind(P^\prime) = -\ind(P)+\ind(P^\prime)=0$$
and the result is proved.
\qed

In the remaining sections we will prove Theorem \ref{cobdsm inv} making use of the so-called 'push-forward map' in $K$-theory, considered by Atiyah and Singer in \cite{a-sI}.

\section{The push-forward map}\label{s p-fwd map}

Let $X$, $Y$ be manifolds, possibly with boundary, and let $i: X\to Y$ be a $K$-oriented embedding with normal bundle $N$ (that is, such that $N$ is even-dimensional and is endowed with a \spinc-structure).
 To define the push-forward map induced by $i$, we require that $X$ have an open tubular neighborhood in $Y$, which will then be identified with the normal bundle $N$. 
Such a neighborhood always exists when $\partial X=\partial Y=\emptyset$; 
 if the boundaries are non-empty, it comes down to 
requiring that the embedding $i$ is \textit{neat}
(see \cite{hirsch}).
 More precisely, we require that

(i) ${X}\cap\partial Y=\partial X$,

(ii) $X$ is not tangent to $\partial Y$ at any point $x\in\partial X$,
that is, $T_x{X}\nsubseteq T_x(\partial Y)$. 

\noindent These conditions are, in fact, equivalent to the existence of an open tubular neighborhood of $X$ in $Y$. 

 Since $N$ is  $K$-oriented, there is a Thom isomorphism given by 
\begin{equation}
\rho : K^0(X) \to K^0(N),\;\;a\mapsto \lambda_N\cup a,
\end{equation}
where $\lambda_N$ is the Thom class (see \cite{atiyah-K, karoubi}). On the other hand, identifying $N$ with an  neighborhood of $X$ in $Y$, we have a $K$-theory map, as in (\ref{ext}),
\begin{equation}
h : K^0(N) \to K^0(Y)
\end{equation}
induced by the open embedding $N\to Y$, that is, by the collapsing map $Y^+\to Y^+/(Y^+-N)\cong N^+$.
\begin{defn}\label{push-fwd}
Let $i: X\to Y$ be a $K$-oriented, neat embedding. The \textup{push-forward map} $i_! : K^0(X) \to K^0(Y)$ is defined as 
$$i_!:= h\circ \rho.$$
\end{defn}
This map does not depend on the open tubular neighborhood chosen \cite{a-sI}.
 We can also define a push-forward map between $K^1$-groups: it is simply the push-forward map associated to the induced embedding $\tilde{i}: X\times\RR\to Y\times\RR$. It is clear that $\tilde{i}$ is $K$-oriented and neat whenever $i$ is, and we then have a push-forward map
\begin{equation}\label{i_!tilde}
\tilde{i}_! : K^1(X) \rightarrow K^1(Y).
\end{equation}
Now note that if $i: X\to Y$ is a neat
embedding, not necessarily $K$-orientable, the induced (neat) embedding of tangent spaces, denoted
again by $i: TX\to TY$ is always $K$-orientable; this follows 
from the fact that the normal bundle of $TX$ in $TY$ is given by $\pi^*(N\oplus N)\cong TN$, where $N$ is the normal bundle of $X$ in $Y$ and $\pi:TX\to X$ is the projection map.
 We have, then that a neat embedding $i:X\to Y$ induces push-forward maps on tangent spaces
\begin{equation}\label{i_! TX}
i_! : K^0(TX) \rightarrow K^0(TY)
\end{equation}
\begin{equation}\label{i_!tilde TX}
\tilde{i}_! : K^1(TX) \rightarrow K^1(TY),
\end{equation}
associated to the  $K$-oriented, neat embeddings $i: TX\to TY$ and $i: TX\times \RR\to TY\times\RR$.
(The map (\ref{i_! TX}) is the one defined by Atiyah and Singer in \cite{a-sI}.) 

The transitivity of the Thom isomorphism yields the functoriality of the push-forward map, that is, if $i: X\to Y$, $j: Y\to Z$ are  $K$-oriented, neat embeddings, then
\begin{equation}\label{funct i_!}
j_!\circ i_! = (j\circ i)_!.
\end{equation}

\noindent We will be interested in analyzing in more detail two particular cases of functoriality:

(i) if $X\subset Y$ is open and $i$ is the inclusion, then $i_!$ coincides with the map induced by an open embedding, as in (\ref{ext});

(ii) if $Y$ is a vector bundle over $X$ with a spin$^c$-structure and $i$ is the zero-section embedding, then $i_!$ coincides with the Thom isomorphism. In particular, if $Y=X\times\RR^2$, then $i_!$ is the Bott isomorphism.


For (i), we note that if an embedding $i :  X \rightarrow Y$ restricts to an embedding $ U \rightarrow W$, with $U\subset X$ and $W\subset Y$ open, that is, if $U=X\cap W$, then the normal bundles to $U$ in $X$ and to $W$ in $Y$ are zero, that is, 
$$TU = TX_{|U},\;\;TW = TY_{|W}.$$
Then 
 the normal bundle to $U$ in $W$ can be identified with $N_{|U}$; in particular, it is $K$-oriented and neat 
 whenever $N$ is. 

\begin{prop}\label{i_! vs ext}
Let $X$, $Y$ be manifolds, $i : X \rightarrow Y$ a  $K$-oriented, neat embedding. Let $U\subset X$ and $W\subset Y$ be open, such that $U=X\cap W$, so that the embedding $i$ restricts to an embedding $i^\prime : U \rightarrow W$. 
Then the following diagram
\[\begin{CD}
K^j(U)   @>{i^\prime_!}>>        K^j(W) \\
 @V{h_X}VV                          @VV{h_Y}V \\
K^j(X)   @>>{i_!}>    K^j(Y),
\end{CD}\]
%
is commutative, for $j=0,1$, where the vertical arrows denote extension maps as in (\ref{ext}).
\end{prop}
\pf
We consider $j=0$ (the case $j=1$ will follow, for $X\times\RR$, $Y\times\RR$, $U\times\RR$, $W\times\RR$). Write $h_X=(j_X)_!$ and $h_Y=(j_Y)_!$, with $j_X: U\to X$, $j_Y: W\to Y$ inclusion maps, and note that $j_Y\circ i^\prime = i\circ j_X$ is just a restatement of $U=X\cap W$. Then (\ref{funct i_!}) yields that $h_Y\circ i^\prime_!=i_!\circ h_X$.
\qed

Noting that $TX\cap TW=T(X\cap W)$, we see that the induced embedding $i:TX\to TY$ also restricts to an embedding $TU\to TW$. Hence, the previous result holds with $X,Y,U,W$ replaced by $TX,TY,TU,TW$ for the maps (\ref{i_! TX}) and (\ref{i_!tilde TX}), where $i$ is not assumed to be $K$-oriented.

We now consider case (ii);  we will consider only the Bott isomorphism. Note that an embedding $i : X \rightarrow Y$ induces an embedding 
 $i^\prime : X \times \RR^2 \rightarrow Y\times \RR^2$  and, it is clear that $i^\prime$ is $K$-oriented and neat whenever $i$ is.

\begin{prop}\label{i_! vs Bott}
Let $X$ and $Y$ be  manifolds and $i : X \rightarrow Y$   a  $K$-oriented, neat embedding. Let  $i^\prime : X \times \RR^2 \rightarrow Y\times \RR^2$ be the induced embedding. 
 Then the following diagram 
\[\begin{CD}
K^j(X)                 @>{i_!}>>              K^0(Y) \\
 @V{\beta_{X}}VV                                @VV{\beta_{Y}}V \\
K^j(X \times \RR^2)   @>>{i^\prime_!}>    K^j(Y\times \RR^2).
\end{CD}\]
is commutative, for $j=0,1$, where  $\beta_X: K^j(X)\to K^j(X\times\RR^2)$, $\beta_Y: K^j(Y)\to K^j(Y\times\RR^2)$ denote the Bott isomorphisms.
\end{prop}
\pf
It is again enough to consider the case $j=0$. Write $\beta_X = (j_X)_!$ and $\beta_X = (j_Y)_!$, where $j_X : X\to X\times\RR^2$ and $j_Y : Y\to Y\times\RR^2$ are the zero-section embeddings. Then it is clear that $j_Y\circ i = i^\prime\circ j_X$ and (\ref{funct i_!}) gives the result.
\qed

Again, the diagram above commutes with $X, Y$ replaced by $TX, TY$, for the maps (\ref{i_! TX}) and (\ref{i_!tilde TX}), where $i$ is not assumed to be $K$-oriented.

We now prove another functoriality result that does not stem directly from (\ref{funct i_!}). 
 For a closed $X_0\subset X$, the inclusion $k : X_0 \rightarrow X$ induces a restriction map
\begin{equation}\label{rest}
r_X := k^* : K^0(X) \rightarrow K^0(X_0). 
\end{equation}
We are particularly interested in restriction to the boundary. Note that a neat embedding $i:X\to Y$ induces an embedding of the boundaries, $i:\partial X\to\partial Y$, which is $K$-orientable whenever $i$ is.

\begin{prop}\label{i_! vs rest to bdry}
Let $X$, $Y$ be manifolds, and $i : X \rightarrow Y$ be a $K$-oriented, neat embedding. The following diagram
\[\begin{CD}
K^j(X)   @>{i_!}>>        K^j(Y) \\
 @V{r_X}VV                          @VV{r_Y}V \\
K^j(\partial X)   @>>{i_!}>    K^j(\partial Y),
\end{CD}\]
is commutative, for $j=0,1$, where the vertical arrows denote restriction maps.
\end{prop}
\pf
We consider only $j=0$ (from which the result for $j=1$ follows). 
Let $N$  denote the normal bundle to $X$ in $Y$; the normal bundle to $\partial X$ in $\partial Y$ is given  by $N_{|\partial X}=\partial N$. Hence we can write the diagram as
\[\begin{CD}
K^0(X)     @>{\rho_{N}}>>   K^0(N)     @>{h_{Y}}>>        K^0(Y) \\
@V{r_X}VV                    @VVV                      @VV{r_Y}V \\
K^0(\partial X)  @>>{\rho_{\partial N}}>  K^0(\partial N)  @>>{h_{\partial Y}}>      K^0(\partial Y),
\end{CD}\]
where the middle vertical arrow is also given by restriction. Functoriality of the Thom isomorphism gives that the left-hand side diagram commutes. As for the right-hand side diagram, note that the composition $r_Y\circ h_Y$ is the map (in reduced $K$-theory) induced by $k\circ i_Y$, where $i_Y: (\partial Y)^+\to Y^+$ is the inclusion and $k: Y^+\to Y^+/Y^+-N =N^+$ is the collapsing map. Since $N$ is open in $Y$, we have that $\partial Y\cap N= \partial N$; therefore $k\circ i_Y$ coincides with $i_N\circ k_0$, where now $i_N: (\partial N)^+\to N^+$ is the inclusion and $k_0: (\partial Y)^+\to (\partial N)^+$ is the collapsing map. This proves the result.
\qed

Concerning the push-forward maps (\ref{i_! TX}) and (\ref{i_!tilde TX}) on tangent spaces,  $i_!:K^0(TX)\to K^0(TY)$ and $\tilde{i_!}:K^1(TX)\to K^1(TY)$, note that since there is an embedding $h:\partial X\times [\,0,1)\to X$ such that $h(x,0)=x$,
 one has that 
\begin{equation*}\label{tg bdle bdry}
TX_{|\partial X}\cong T(\partial X) \times \RR
\end{equation*}
(where we take the inward normal vector at each point of $\partial X$). Under this isomorphism, one has that the push-forward map $i_!:K^0(TX_{|\partial X})\to K^0(TY_{|\partial Y})$ gives the map (\ref{i_!tilde TX}), $\tilde{i_!}:K^1(T(\partial X))\to K^1(T(\partial Y))$. This remark, together with Proposition \ref{i_! vs rest to bdry}, applied to $TX$, $TY$, gives the following.

\begin{cor}\label{i_! TX vs rest}
Let $X$, $Y$ be  manifolds and $i : X \rightarrow Y$ a neat embedding. 
Then the following diagram
\[\begin{CD}
K^j(TX)   @>{i_!}>>        K^j(TY) \\
 @V{r_X}VV                          @VV{r_Y}V \\
K^{j+1}(T(\partial X))   @>>{\tilde{i_!}}>    K^{j+1}(T(\partial Y)),
\end{CD}\]
is commutative, for $j = 0,1$, where the vertical arrows denote restriction maps composed with the isomorphisms $K^j(TX_{|\partial X})\cong K^j(T(\partial X)\times\RR)$.
\qed
\end{cor}

One of the crucial features of the push-forward map (\ref{i_! TX}) in the index theory setting is its compatibility with the analytical index map on compact manifolds without boundary. The (highly non-trivial) proof of the following result can be found in \cite{a-sI}.

\begin{thm}[Atiyah-Singer]\label{i_! vs an-ind cpt}
Let $X$, $Y$ be closed manifolds, and $i: X \rightarrow Y$ an embedding. Let $a \in K^0(TX)$ and consider $i_! : K^0(TX) \rightarrow K^0(TY)$. Then $$\ind (i_!(a)) = \ind (a).$$
\end{thm}

We will now show that the same is true in the non-compact case, where we will be computing the indices of  multiplication operators at infinity.  
 We first prove this result in a nice particular case.  In the remaining of this section, $M$ and $W$ are $\sigma$-compact manifolds without boundary, and $i:M\to W$ is an embedding (necessarily neat). All our embeddings have closed image.

\begin{lem}\label{anind vs i_! emb in cpt w bdry}
 Let $M$, $W$ be manifolds, $i:M\to W$ an embedding and consider the push-forward map $i_! : K^0(TM) \rightarrow K^0(TW)$. If there exist compact manifolds  $X$ and $Y$, possibly with boundary, such that $M\subset X$, $W\subset Y$ are open, and $i$ extends to a neat embedding $i^\prime:X\to Y$, then $\ind \circ i_! = \ind$.
\end{lem}
\pf
 Assume first that $\partial X=\partial Y=\emptyset$. Consider the  maps $$h_X: K^0(TM)\to K^0(TX),\;h_Y: K^0(TW)\to K^0(TY)$$ associated to the respective open inclusions. For any $a\in K^0(TM)$, using excision (Proposition \ref{an-ind vs ext}) and Proposition \ref{i_! vs ext} for tangent spaces, we have that
$$\ind(i_!(a))=\ind(h_{Y}(i_!(a)))=\ind(i^\prime_!(h_{X}(a))).$$
Applying Theorem \ref{i_! vs an-ind cpt} to $X$, $Y$, which are compact without boundary, we get
$$\ind(i^\prime_!(h_{X}(a)))=\ind(h_{X}(a))=\ind(a),$$
where we have again used excision.

When $X$ and $Y$ have a boundary, let $X^d$, $Y^d$ be the doubles\footnote{If $\partial X\neq\emptyset$, we define the double of $X$ as the smooth manifold $X^d$ obtained by taking the union of two copies of $X$ and identifying their boundaries.} of $X$, $Y$, respectively; $X^d$, $Y^d$ are compact manifolds without boundary. Since the embedding $i^\prime:X\to Y$ is neat, we have that $\partial Y\cap X=\partial X$ and we can extend $i^\prime:X\to Y$, and therefore also $i:M\to W$, to an embedding $i^d: X^d\to Y^d$. 
 Moreover, $X\subset X^d$ and $Y\subset Y^d$ are open, and hence $M\subset X^d$ and $W\subset Y^d$ are open as well. We conclude that the result holds for any compact manifolds $X$, $Y$.
\qed

It is not true in general that an open manifold can be embedded in a compact manifold (with or without boundary). However, using Sard's theorem, one can show that any $\sigma$-compact manifold can be covered by open submanifolds with this property.
  
\begin{lem}\label{lem sigmacpt}
Any $\sigma$-compact manifold $M$ can be written as $M = \bigcup_{n\in\NN} M_n$ where $M_n$ is open, $\overline{M}_n\subset M_ {n+1}$ and $\overline{M}_n$ is a compact submanifold of $M$ (possibly with boundary). 
\qed
\end{lem}
In this case, we have that
\begin{equation}\label{dirlim-}
K^0(TM) = \lim\limits_{\to} K^0(TM_n),
\end{equation}
with morphisms $h_{n} : K^0(TM_n) \to K^0(TM)$,  given by the usual maps associated to open inclusions, as in (\ref{ext TX}). One more application of excision will then enable us to reduce the computation of the index from $M$ to some $M_n$.

\begin{thm}\label{i_! vs an-ind}
Let $M$, $W$ be manifolds, $i: M \rightarrow W$ be an embedding and consider the push-forward map $i_! : K^0(TM) \rightarrow K^0(TW)$. Then, $$\ind \circ i_! = \ind.$$
\end{thm}
\pf
Write $i_! = h\circ\rho$, where $\rho$ is the Thom isomorphism and  $h$ is the $K$-theory map associated to an open inclusion. We have shown in Proposition \ref{an-ind vs ext} that $\ind\circ h=\ind$.
 We can therefore assume that $i=\rho$, that is, that $W$ is a vector bundle over $M$ and $i: M\to W$ is the zero-section embedding; we let $\pi : W \to M$ be the projection map.
 
Write $M=\bigcup_{n\in\NN}M_n$, as in Lemma \ref{lem sigmacpt}, and take $a_n\in K^0(TM_n)$ such that $a=h_n(a_n)$. Defining $W_n := \pi^{-1}(M_n)$, we have that $W_n$ is a vector bundle over $M_n$ and that $i$ restricts to the zero-section embedding $i_n : M_n \to W_n$. Using excision, Proposition \ref{an-ind vs ext}, and the functoriality of the push-forward map, Proposition \ref{i_! vs ext}, we have that
$$\ind(a)=\ind(a_n)\quad\text{and}\quad\ind(i_!(a))=\ind({i_n}_!(a_n)),$$
and hence, to prove our claim, it suffices to consider the embedding $i_n: M_n \to W_n$, that is, we can assume without loss of generality that $M$ is open such that $\overline{M}$ is a compact manifold (with boundary). 

Now, $i:M\to W$ clearly extends to an embedding $i:\overline{M}\to\overline{W}$, and, moreover, $M\subset \overline{M}$,  $W\subset \overline{W}$ are open, with $\overline{M}_n$ compact. (Note that  $i$ is neat, since $\overline{W}$ is a vector bundle over $\overline{M}$.) Now define a fiber bundle $E$ over $\overline{M}$ such that $E_x= ({\overline{W}})_x^+$, where $({\overline{W}})_x^+$ is the one-point compactification of the fiber $({\overline{W}})_x$, $x\in \overline{M}$. We conclude that we can extend $i$ to a neat embedding 
$i^\prime: \overline{M}\to E$, where $\overline{M}$, $E$ are compact manifolds with boundary, and $M\subset\overline{M}$, $W\subset E$ are open.  Lemma \ref{anind vs i_! emb in cpt w bdry} then yields the result in full generality.
\qed

\section{Proof of cobordism invariance}\label{s pf cbd inv}

We now prove Theorem \ref{cobdsm inv}. We start with showing that cobordism of symbols is preserved by the push-forward map associated to an embedding. Using Theorem \ref{i_! vs an-ind}, we reduce the proof to $\RR^n$ where it is seen to be trivial. We first analyze the behavior of the push-forward map with respect to the symbol restriction map (\ref{u_M}). Recall that $\tilde{i_!}:K^1(TX)\to K^1(TY)$, as in (\ref{i_!tilde TX}), is the push-forward map induced by the embedding $TX\times\RR\to TY\times\RR$.

\begin{lem}\label{i_! vs u_M}
Let $X$, $Y$ be manifolds, $i : X \rightarrow Y$ a neat embedding and denote the boundaries of $X$ and $Y$ by  $X_0$ and $Y_0$, respectively. 
Let $$u_{X_0} = \beta_{TX_0}^{-1} \circ r_{X_0},\;\text{and}\; u_{Y_0} = \beta_{TY_0}^{-1} \circ r_{Y_0}$$ be as in (\ref{u_M}). Then the following diagram is commutative:
\[\begin{CD}
K^1(TX)   @>{\widetilde{i_!}}>>        K^1(TY) \\
 @V{u_{X_0}}VV                          @VV{u_{Y_0}}V \\
K^0(TX_0)   @>>{i_!}>    K^0(TY_0).
\end{CD}\]
\end{lem}
\pf
From Corollary \ref{i_! TX vs rest}, we have that
\[\begin{CD}
K^1(TX)   @>{\widetilde{i_!}}>>        K^1(TY) \\
@V{r_{X}}VV                                             @VV{r_{Y_0}}V \\
K^1(TX_0\times \RR)   @>>{\widetilde{i_!}}>        K^1(TY_0\times \RR)
\end{CD}\]
commutes. On the other hand, since $i$ is neat, it induces a smooth embedding $i : X_0 \rightarrow Y_0$ and, from Proposition \ref{i_! vs Bott} for tangent spaces, the diagram
\[\begin{CD}
K^0(TX_0\times \RR^2)   @>{i_!}>>   K^0(TY_0 \times \RR^2) \\
@V{{\beta}^{-1}_{TX_0}}VV            @VV{{\beta}^{-1}_{TY_0}}V \\        
K^0(TX_0)                @>>{i_!}>           K^0(TY_0)
\end{CD}\]
commutes as well. The map $i_! : K^0(TX_0\times \RR^2)\to K^0(TY_0 \times \RR^2)$  in the top row is the push-forward map induced by the embedding $i :
X_0 \times \RR \rightarrow Y_0\times \RR$, which is easily seen to coincide with the definition of $ {\widetilde{i_!}}$.
We can then couple the diagrams together to get the result.
\qed

\begin{prop}\label{i_! vs cobdsm symb}
Let $M$ be a manifold and let $\sigma\in K^0(TM)$ be such that $(M,\sigma)\sim 0$ for a cobordism of symbols $(X,\omega)$. Let $Y$ be a manifold, and $i:X\to Y$ a neat embedding. Then, if  $M^\prime=\partial Y$, we have $(M^\prime, i_!(\sigma))\sim 0$ with the cobordism of symbols given by $(Y, \,\tilde{i}_!(\omega))$.
\end{prop}
\pf
From the previous lemma, with $X_0=M$, $Y_0=M^\prime$, if $\sigma\in K^0(TM)$, $\omega \in K^1(TX)$ are such that $u_M(\omega)=\sigma$, then $u_{M^\prime}(\tilde{i}_!(\omega))=i_!(\sigma)$.
\qed

We will use this result for an embedding of $X$ in a manifold $Y = [0,+\infty
) \times \RR ^n$  such that $\partial X=M$ is embedded in $\{0\} \times \RR ^n$. Using the invariance of the analytical index with respect to the push-forward map (Theorem \ref{i_! vs an-ind}) we are then able to reduce the proof of cobordism invariance to operators defined on $\RR^n$. In this particular case, the proof is straightforward due to our second lemma.

\begin{lem}\label{K=0}
Let $Z = [\,0,+\infty) \times \RR ^n$. Then $K^1(TZ) = 0$. 
\end{lem}
\pf
We have $TZ = [\,0,+\infty)\times\RR^{2n+1}$ and therefore $(TZ\times\RR)^+\cong B^{2n+2}$. Since $B^{2n+2}$ is contractible, $K^1(TZ)=\tilde{K^0}(B^{2n+2})=0$.
\qed

We can now finally prove that the Fredholm index of a bounded elliptic \pso that is multiplication at infinity is invariant under cobordism of symbols.
\newline

\noindent\textbf{Proof of Theorem \ref{cobdsm inv}.}
Let $P\in\Psi_{mult}^b(M;E,F)$ be a bounded elliptic operator and let $\sigma:=[\sigma_0(P)] \in K^0(TM)$ denote the principal symbol class of $P$. Assuming that $(M,\sigma)\sim 0$, we want to show that $\ind(\sigma)=\dim\ker(P)-\dim\coker(P) = 0.$

Let  then $X$ be a  manifold with $M=\partial X$ and $\omega \in K^1(TX)$ such that $u_M(\omega)=\sigma$.
Consider an  embedding $i : X \rightarrow Y$ with  $Y = [0, + \infty[ \times
\RR^n$ such that $\partial X=M$ is embedded in $\{0\} \times \RR ^n$. In this case, we can always assume that $i$ is neat (see \cite{hirsch}, Theorem 1.4.3, for the proof in the compact case; the non-compact case can be proved in exactly the same way). From Proposition \ref{i_! vs cobdsm symb}, with $M^\prime= \RR^n$,  we have that $(M^\prime, i_!(\sigma))\sim 0$, and moreover, that  
$$i_!(\sigma) = u_{M^\prime}(\widetilde{i_!}(\omega)).$$ 
 By Lemma \ref{K=0}, we know that $\tilde{i}_!(\omega) = 0$ and therefore $i_!(\sigma) = 0$. On the other hand, from 
 Theorem \ref{i_! vs an-ind}, we conclude that
$$\ind(\sigma) = \ind(i_!(\sigma)) = 0.$$
\qed

\section{An index formula}\label{s ind form}

We show here how Theorem \ref{i_! vs an-ind} can be used to extend the $K$-theoretical index formula given in \cite{a-sI} to \psos that are multiplication at infinity. 
Let $M$ be a $\sigma$-compact manifold without boundary and consider an embedding $i : M \to \RR^n$, for some $n\in\NN$ (such an embedding always exists, see \cite{hirsch}). Let $j : P_0 \to \RR^n$ be the inclusion of a point; then $j_! : K^0(TP_0) = \ZZ \to K^0(T\RR^n)$ coincides with the Bott isomorphism. 
 The \textup{topological index} map $\topind: K^0(TM) \to  \ZZ$ is defined as
\begin{equation*}
\topind := (j^{-1}_!)\circ i_!.
\end{equation*}

\begin{thm}\label{ind form}
Let $P\in\Psi_{mult}^0(M;E,F)$ be an bounded elliptic operator and $\sigma:=[\sigma_0(P)] \in K^0(TM)$ denote its principal symbol class. Then $$\ind(P) = \topind(\sigma).$$
\end{thm}
\pf
Let $\ind^M$, $\ind^{\RR^n}$ and $\ind^{pt}$ denote the index maps on $M$, $\RR^n$ and $P_0$, respectively. It is clear that $\ind^{pt}=id$.
 Write $\topind^M= (j_!)^{-1}\circ i_!$, where $i: M\to \RR^n$ is an embedding and $j: P_0\to \RR^n$ is the inclusion of a point. From Theorem \ref{i_! vs an-ind}, the diagram
\[\begin{CD}
K^0(TM)         @>{i_!}>>   K^0(T\RR^n)        @<{j_!}<<  K^0(TP_0)=\ZZ  \\
@V{\ind^{M}}VV            @VV{\ind^{\RR^n}}V           @VV{\ind^{pt}}V\\
\ZZ             @=           \ZZ            @=         \ZZ
\end{CD}\]
commutes. Since $\ind^{pt} = id$, we conclude that $\ind^{M} = \topind^M$.
\qed

\noindent Note that, from the proof of this result, we have that the index is totally characterized by
\newline
\indent(i) $\ind^{pt}=id$,
\newline
\indent(ii) $\ind\circ i_!=\ind$,
\newline
as it is on closed manifolds \cite{a-sI}.
\newline


\end{document}